\documentclass[reqno,11pt]{amsart}
\usepackage{amsmath, latexsym, amsfonts, amssymb, amsthm, amscd}
\usepackage{psfig}
\usepackage{graphics,epsf,psfrag}
\numberwithin{equation}{section}

\newtheorem{theorem}{Theorem}[section]
\newtheorem{lemma}[theorem]{Lemma}
\newtheorem{proposition}[theorem]{Proposition}

\newtheorem{rem}[theorem]{Remark}
\newtheorem{definition}[theorem]{Definition}

\DeclareMathOperator{\sign}{\mathrm{sign}}

\newcommand{\ind}{\mathbf{1}}

\newcommand{\R}{\mathbb{R}}

\newcommand{\N}{\mathbb{N}}
\renewcommand{\tilde}{\widetilde}

\newcommand{\cA}{{\ensuremath{\mathcal A}} }

\newcommand{\cN}{{\ensuremath{\mathcal N}} }
\newcommand{\cL}{{\ensuremath{\mathcal L}} }

\newcommand{\cD}{{\ensuremath{\mathcal D}} }

\newcommand{\bP}{{\ensuremath{\mathbf P}} }
\newcommand{\bE}{{\ensuremath{\mathbf E}} }

\DeclareMathSymbol{\leqslant}{\mathalpha}{AMSa}{"36} 
\DeclareMathSymbol{\geqslant}{\mathalpha}{AMSa}{"3E} 
\DeclareMathSymbol{\eset}{\mathalpha}{AMSb}{"3F}     
\newcommand{\dd}{\,\text{\rm d}}             

\newcommand{\bbE}{{\ensuremath{\mathbb E}} }

\newcommand{\bbL}{{\ensuremath{\mathbb L}} }

\newcommand{\bbP}{{\ensuremath{\mathbb P}} }

\newcommand{\gb}{\beta}

\newcommand{\gep}{\varepsilon}       

\newcommand{\gD}{\Delta}

\newcommand{\go}{\omega}
\newcommand{\gO}{\Omega}
\newcommand{\gl}{\lambda}

\makeatletter
\def\captionfont@{\footnotesize}
\def\captionheadfont@{\scshape}

\long\def\@makecaption#1#2{%
  \vspace{2mm}
  \setbox\@tempboxa\vbox{\color@setgroup
    \advance\hsize-6pc\noindent
    \captionfont@\captionheadfont@#1\@xp\@ifnotempty\@xp
        {\@cdr#2\@nil}{.\captionfont@\upshape\enspace#2}%
    \unskip\kern-6pc\par
    \global\setbox\@ne\lastbox\color@endgroup}%
  \ifhbox\@ne 
    \setbox\@ne\hbox{\unhbox\@ne\unskip\unskip\unpenalty\unkern}%
  \fi
  \ifdim\wd\@tempboxa=\z@ 
    \setbox\@ne\hbox to\columnwidth{\hss\kern-6pc\box\@ne\hss}%
  \else 
    \setbox\@ne\vbox{\unvbox\@tempboxa\parskip\z@skip
        \noindent\unhbox\@ne\advance\hsize-6pc\par}%
\fi
  \ifnum\@tempcnta<64 
    \addvspace\abovecaptionskip
    \moveright 3pc\box\@ne
  \else 
    \moveright 3pc\box\@ne
    \nobreak
    \vskip\belowcaptionskip
  \fi
\relax
}
\makeatother
\def\writefig#1 #2 #3 {\rlap{\kern #1 truecm
\raise #2 truecm \hbox{#3}}}

\newcommand{\tf}{\textsc{f}}
\newcommand{\M}{\textsc{M}}

\newcommand{\rc}{\mathtt c}
\newcommand{\rf}{\mathtt f}
\newcommand{\Ndown}{\cN}
\newcommand{\overh}{\overline{h}}

\begin{document}

\title[Estimates on path delocalization ]{
Estimates on path delocalization \\ for copolymers at selective interfaces 
}
\author{Giambattista Giacomin}
\address{Laboratoire de Probabilit{\'e}s de P 6\ \& 7 (CNRS U.M.R. 7599)
  and  Universit{\'e} Paris 7 -- Denis Diderot,
U.F.R.                Mathematiques, Case 7012,
                2 place Jussieu 75251 Paris cedex 05, France
\hfill\break
\phantom{br.}{\it Home page:}
{\tt http://felix.proba.jussieu.fr/pageperso/giacomin/GBpage.html}}
\email{giacomin\@@math.jussieu.fr}
\author{Fabio Lucio Toninelli}
\address{Institut f\"ur Mathematik, Universit\"at Z\"urich,
Winterthurerstrasse 190, CH--8057 Z\"urich, Switzerland
\hfill\break
Present address: Laboratoire de Physique, UMR-CNRS 5672, ENS Lyon, 46 All\'ee d'Italie, 69364 Lyon Cedex 07, France}
\email{fltonine@ens-lyon.fr}
\date{\today}

\begin{abstract}
Starting from the simple symmetric random walk $\{ S_n \}_n$, we introduce 
a new process whose path measure is weighted by a factor
$\exp\left( \gl \sum_{n=1}^N \left(\go_n +h \right ) \sign \left( S_n\right)\right)$,
with $\gl , h \ge 0$,  $\{ \go _n \}_n $ a typical realization
of an IID process and $N$ a positive integer. We are looking for results in the large $N$ limit.
This factor favors $S_n>0$ if $\go_n+h >0$
and $S_n<0$ if $\go_n+h <0$. 
The process can be interpreted
as a model for a random
heterogeneous polymer in the proximity of an
interface separating two selective solvents.
It has been shown  \cite{cf:BdH} that
this model undergoes a (de)localization transition:
more precisely there exists a continuous increasing function
$\gl \longmapsto h_c(\gl)$ such that
if $h< h_c(\gl)$ then the model is localized while it is delocalized if $h\ge h_c(\gl)$.
However, localization and delocalization were not given in terms
of path properties, but in a free energy sense.
Later on it has been shown that free energy localization
does indeed correspond to a (strong) form of path localization \cite{cf:BisdH}.
On the other hand, only weak results on the delocalized regime have been known so far.
\\
We present a method, based on concentration bounds on
{\sl suitably restricted} partition functions, that yields much stronger 
results on the path behavior in the interior of the delocalized region,
that is for $h> h_c (\gl)$. In particular we prove that, in a suitable sense,   
one cannot expect more than $O(\log N)$ 
visits of the walk to the lower half plane. The previously known bound was $o(N)$.
Stronger $O(1)$--type results are obtained deep inside the delocalized region.

The same approach is also helpful for a different type of question: we prove in fact that
the limit as $\gl$ tends to zero of $h_c(\gl) / \gl$ exists and it is
independent of the law of $\go _1$, at least when the random variable  $\go_1$
is bounded or it is Gaussian. This is achieved by  interpolating between 
this  class of variables and the particular case of $\go_1$ taking
values $\pm 1$ with probability $1/2$, treated in \cite{cf:BdH}.
\end{abstract}

\maketitle

 \section{Introduction}
 \label{sec:intro}
\setcounter{equation}{0}

\subsection{The model and its free energy}
Let $S=\{S_n \}_{n=0,1,\ldots}$ be a simple random walk:
 $S_0=0$ and  $\{ S_j- S_{j-1}\}_{j\in\N}$ a sequence
of IID random variables with  $\bP \left( S_1= \pm 1\right)=1/2$. 
We denote by $\gO$ the set of all random walk trajectories.
For $\gl\ge 0$,  $h \ge 0$, $N\in 2\N$  and 
$\go =\{ \go_n\}_{n=1,2, \ldots} \in \R ^\N$
we introduce the {\sl copolymer} measures 
\begin{equation}
\label{eq:Boltzmann}
\frac {\dd \bP_{N, \go}^{a}} {\dd \bP} (S)
\, = \, 
\frac 1
{\tilde Z_{N,\go}^{a}}
{\exp\left( \gl \sum_{n=1}^N \left( \go_n +h\right) \sign \left(S_n\right)
\right)} \ind_{\gO_N^a}
,
\end{equation}
with $a=\rf$ (free case) or $a=\rc$ (constrained case),
$\gO_N^\rf=\gO$, 
$\gO_N ^\rc = \{S\in \gO:  S_N=0\}$.
  $\tilde Z_{N,\go}^{a}$ is the partition function and 
$\sign \left(S_{2n}\right)$ is set to be equal to $\sign \left(S_{2n-1}\right)$
for any $n$ such that $S_{2n}=0$. 

\smallskip

The sequence $\go$ is chosen as a typical realization of an IID
sequence of random variables, still denoted by $\go = \{ \go_n \}_n$.
We call $\bbP$ the law of $\go$.
Further hypotheses on $\go$ are summed up by: 
\smallskip

\begin{definition}
\label{def:omega}
$\empty$
\begin{itemize}
\item {\bf Basic assumptions}:   $\go_1 \sim -\go_1$ and $\M (t):= 
\bbE\left[ \exp\left(t \go_1\right)\right]< \infty$ for $ t$ in a neighborhood of zero.
Without loss of generality we assume $\bbE \left[ {\go_1} ^2 \right]=1$.
\item {\bf Deviation inequality above the mean}: there exists
 a positive constant $C$ such that for every $N$,
 for every Lipschitz and convex function $g: \R ^ N \to \R$ with
 $g(\go) := g\left(\go_1, \ldots, \go _N \right) \in {\mathbb L}^1 \left( \bbP \right)$ and  $t\ge 0$
 \begin{equation}
 \label{eq:Lip}
 \bbP 
 \left(  g(\go) - \bbE \left[ g(\go)\right] \ge t\right)
 \, \le \,  C \exp \left( -\frac{t^2}{C \Vert g \Vert ^2 _{\rm Lip}}\right),
 \end{equation}
 where
  $\Vert g \Vert_{\rm Lip}$ is the Lipschitz constant of $g$ with respect to the Euclidean distance.
\end{itemize}
\end{definition}
\smallskip

The deviation inequality (\ref{eq:Lip}) is known to hold with
a certain generality: its validity for
 the Gaussian case $\omega_1\sim{\mathcal N(0,1)}$
and for the case of bounded random variables is by now 
a classical result, see \cite{cf:newlook}, 
\cite{cf:Ledoux} and \cite{cf:cedric}. However one can go beyond: it holds in particular
whenever the law of $\go_1$ satisfies the log--Sobolev inequality
\cite{cf:Ledoux} and in that case of course $C$ depends on the log--Sobolev constant. 
As a matter of fact, in all the cases we have mentioned not only
a deviation inequality above the mean holds, but also below, and therefore
one has the full concentration inequality. 
A necessary and sufficient condition for 
the log--Sobolev inequality to hold can be found in 
\cite{cf:BoGoe}. 
In order  to be more explicit we point out that  if $\go_1$
has a density of the type $\exp(-V)$, with $V$ bounded from below 
and strictly convex outside a finite interval, the law of $\go_1 $
satisfies the log--Sobolev inequality with a finite constant
 and therefore \eqref{eq:Lip} holds.

\smallskip

Under the basic assumptions on $\go$
the {\sl quenched free energy} of the system exists, namely
the limit
\begin{equation}
\label{eq:free_energy}
f(\gl ,h)\, := \, \lim_{N\to \infty} \frac 1N \log \tilde Z^{a}_{N,\go},
\end{equation}
exists in the $\bbP\left( \dd \go \right)$--almost
sure sense and in the $\bbL^1\left( \bbP\right)$ sense.  
This existence result 
can be proven via super--additivity arguments (we refer to \cite{cf:G} for the details) and the method
shows also that $f(\gl,h)$ is non-random and independent
of the choice of $a$. 
 
We observe that
\begin{equation}
\label{eq:delocfe}
f(\gl , h) \, \ge \, \gl h.
\end{equation} 
The proof of such a result is elementary:
 if we set $\gO_N^+= \{S\in \Omega:\, S_n>0$ for $ n=1, 2, \ldots , N\}$ 
 we have
\begin{multline}
\label{eq:step_deloc}
\frac 1N \log \tilde Z_{N,\go}^{\rf}  \ge  
\frac 1N \log \bE 
\left[ 
\exp\left(
\gl \sum_{n=1}^N 
\left( \go_n +h\right) \sign \left(S_n\right)
\right)
;  \gO_N^+
\right]
\\
= \frac {\gl} {N} \sum_{n=1}^N \left( \go_n +h\right) \, + \, 
\frac 1N 
\log \bP \left( \gO_N^+\right)\, \stackrel{N \to \infty}{\longrightarrow}\, \gl h ,
\end{multline}
where the limit is taken  in the almost sure sense:
we have  applied the strong law of large numbers and  the
well known fact that $\bP\left( \gO_N^+\right)$ behaves like 
$N^{-1/2}$ for $N$ large \cite[Ch. 3]{cf:Feller1}. 
The observation \eqref{eq:delocfe}, above all if viewed in the light
of its proof, suggests the following
partition of the parameter space (or {\sl phase diagram}):

\smallskip
\begin{itemize}
\item The localized region: $\cL = \left\{ (\gl , h): \, f(\gl, h)> \gl h\right\}$;
\item The delocalized region:  $\cD = \left\{ (\gl , h): \, f(\gl, h)= \gl h\right\}$.
\end{itemize}
\smallskip

We sum up the known results on the phase diagram:

\medskip

\begin{theorem}
\label{th:main_fe} 
Under the basic assumptions on $\go$
there exists  an increasing function $h_c: [0, \infty )
\longrightarrow  [0, \infty]$ such that  
\begin{equation}
\cL = \left\{ (\gl, h): \, h< h_c(\gl )\right\}
\ \ \text{ and } \ \ \ 
\cD = \left\{ (\gl, h): \, h\ge  h_c(\gl )\right\}.
\end{equation}
$h_c(\cdot)$ is continuous if it takes values in $[0, \infty)$,
otherwise it is continuous in $[0, \sup\{\gl: h_c(\gl)<\infty\})$.
 Moreover
\begin{equation}
\label{eq:main}
\underline h (\gl) \, :=\, 
\frac 1{4\gl /3} \log \M\left( 4\gl /3\right)\le
h_c (\gl)\le \frac 1{2\gl } \log \M\left( 2\gl \right)
\, =: \, \overline{h} (\gl).
\end{equation} 
\end{theorem}

\medskip

Part of the results in 
Theorem \ref{th:main_fe}  have been  proven in \cite{cf:BdH}. The present
version takes into account the improvements brought by 
\cite{cf:BG}. For the rest of the paper we will refer to $\{(\gl,h): \, h>\overh(\gl)\} \subset \cD$
as {\sl strongly delocalized region}.

The bounds in \eqref{eq:main}
yield that $2/3 \le \liminf_{\gl \searrow 0} h_c(\gl)/\gl  $ and
$ \limsup_{\gl \searrow 0} h_c(\gl)/\gl \le 1 $. In \cite{cf:BdH}
it has been shown that the limit of $h_c(\gl)/ \gl$
exists in the particular case of $\go_1$ taking values $\pm 1$
and it can be expressed in terms of a suitable Brownian copolymer,
suggesting thus a universality of this result. The techniques 
we develop allow to interpolate between the $\pm 1$ case
and more general cases, namely:
\medskip

\begin{theorem}
\label{univers}
The slope of the critical curve at the origin,
\begin{eqnarray}
  m_c\, :=\, \lim_{\lambda\searrow 0}\frac{h_c(\lambda)}\lambda,
\end{eqnarray}
exists and does not depend on the law of $\omega_1$, provided that
$\go_1$ is either a bounded symmetric variable of unit variance
or a standard Gaussian variable.
\end{theorem}

\medskip

Theorem \ref{univers}
 is proven in Section \ref{sec:exapp}.
It is  in the line of the interpolation results
\cite{cf:limite} and \cite{cf:ch}, but here one needs to have a more
explicit control of the $\gl$ dependence  of the error made in the interpolation
procedure. It turns out that the approach that we propose here 
for path estimates yields also this control.
It would be interesting to investigate whether a suitable refinement of the strategy we propose
or an extension of the approach in \cite{cf:BdH}, or possibly
a combination of both, would allow to 
obtain a better result, removing the rather unnatural boundedness requirement on the 
random variables, which arises from our application of the interpolation method.

\subsection{From free energy to path behavior }
The polymer measures $\bP_{N, \go}^a$ have been introduced
in \cite{cf:Sinai} and \cite{cf:BdH} motivated by earlier 
theoretical physics works, in particular  by \cite{cf:GHLO} (for updated
physics developments see \cite{cf:Monthus} and references therein).
It is a model for an heterogeneous polymer, constituted
by charged units ({\sl monomers}). The polymer lives in
a solvent which is also heterogeneous: it is made  
of two solvents in a  state
in which a flat interface is present (an example
familiar to everybody is the case of an oil/water interface).
The sign of the charge determines the preference of
a monomer for one solvent or the other and the absolute value
of the charge plays a role in the intensity of such a preference.
Moreover, in general the situation may be asymmetric:
there may be more charges of a certain sign or
the intensity of the  solvent--monomer interaction 
may not be invariant under the change of sign of the charge
(we are modeling this second situation and $h$ is the asymmetry parameter). 
What we want to analyze is which of the two following scenarios prevails:

\smallskip
\begin{enumerate}
\item The polymer places {\sl most of} the monomers in their preferred
solvent (in the model the $n^{\text{th}}$--monomer is preferably above
the $x$--axis, that plays the role of the interface, if $\go_n+h>0$,  and
below if $\go_n +h<0$). This forces of course the polymer to stick
close to the interface and this is the intuitive concept of a localized polymer path.
\item  The polymer lies almost fully in one of the two solvents. Intuitively
that may happen in an asymmetric case. In such a situation one expects
the polymer to wander away from the interface, since it would be undergoing
a repulsion effect of {\sl entropic origin}: the trajectories staying close
to the interface are very few with respect to the trajectories exploring 
freely a half--space. This is for us a delocalized behavior. 
\end{enumerate}
\smallskip

In principle there is a third reasonable scenario: the case in which
the polymer has large fluctuations between the two solvents. It turns out
that, at least if we disregard the critical case $h=h_c(\gl)$, 
 this situation is possible only in the trivial $\gl =0$ case. 
Moreover scenario (1) is effectively observed if $(\gl,h)\in \cL$
and scenario (2) is verified at least in the interior of $\cD$.
But let us be more precise and let us sum up the state of the art on this issue:

\smallskip
\begin{enumerate}
\item If $(\gl, h) \in \cL$ then very strong localization results are available.
The keyword in this case is {\sl tightness} and one should really think 
of a path essentially as being at distance $O(1)$ from the interface.
  The precise statements are rather involved, due to the presence
  of atypical finite stretches in any typical $\go$, and we prefer to refer to
  \cite{cf:Sinai}, \cite{cf:AZ} and  \cite{cf:BisdH}.
\item The results in the delocalized regime are much more meager. 
All the same the following result is available \cite{cf:BisdH}:
if $(\gl, h) \in \stackrel{\circ}{\cD}$ then
for every $L$
\begin{equation}
\label{eq:delocw}
\lim_{N \to \infty}
\frac 1N \sum_{n=1}^N \bP_{N, \go}^a \left( S_n>L\right) =1, \ \ \ \ \bbP (\dd \go)-\text{a.s.}. 
\end{equation} 
\end{enumerate}
\smallskip

While the delocalization result \eqref{eq:delocw}
is in sharp contrast with the localized scenario (1), it is far from matching
the very strong delocalization results available in polymer models without disorder
(for example in the well known $(1+1)$--dimensional wetting models,
see, e.g., \cite{cf:G}, \cite{cf:IY} and \cite{cf:DGZ}): loosely stated one expects that a typical
delocalized path in the limit of $N\to \infty$ has only a finite number of visits
to the lower half--plane and, as a consequence, a Brownian scaling
result should hold with convergence to well known processes like
the Brownian meander or the $\text{Bessel}(3)$ bridge according to whether $a=\rf$
or $a=\rc$ (see Section \ref{sec:further} for more precision on this issue).   
These are reasonable conjectures, supported also
 by the fact that
in the localized regime the results in the disordered model
match what one observes in the non-disordered case.
One should however stress the essential difference between
the localized and delocalized regions: in the first case
one is in a large deviation regime -- $\bP_{N, \go}^a$
charges a set of trajectories which has exponentially
small probability with respect to $\bP$ -- while
this is not the case in the delocalized regime.
The large deviation machinery does not seem to go 
beyond results of the type \eqref{eq:delocw}.

\smallskip

The purpose of this paper is to present an approach, based on concentration
inequalities, that yields results that go well beyond the {\sl density} result 
  \eqref{eq:delocw}.
  
  \medskip
  
  In order to state our main theorem we need some notation:
  we set $\Delta_n = (1-\sign( S_n))/2$ and introduce the random variable 
  $\Ndown = \sum_{n=1}^N \Delta_n$, counting how many monomers
  are in the lower half--plane (unfavorable solvent).
We introduce also the random set $\cA :=
\{ n\le N: \Delta_n =1 \}\cup \{0\}$ and note that  the (even) number $\max \cA$
identifies the point of last exit of $S$ from the lower half--plane. 

\medskip

\begin{theorem}
\label{th:main}
Under the basic assumptions 
on $\go$ we have that 
\begin{enumerate}
\item if $h > \overh(\gl)$ 
there exists $c$ such that
\begin{equation}
\label{eq:withoutconc1}
\bbE\,
\bP^\rf_{N,\go} \left( \max \cA \le \ell   \right)
\ge 1-  c /\sqrt{\ell+1},
\end{equation}
for every $N$ and every non-negative integer $\ell \le N$. Analogously,
\begin{multline}
\label{eq:withoutconc2}
\phantom{movemov}\bbE\,
\bP ^\rc _{N,\go }
\left( \max \{\cA \cap [0 , N/2]\} \le \ell _1 \text{or } 
\min \{\cA \cap [N/2, N]\} \ge N- \ell_2  
\right) \, \ge
\\ 1- \frac c{\sqrt{\ell_1 \ell_2+1}},
\end{multline}
for every $N $, every $\ell_1 \le N/2$ and 
$\ell_2 \le N/2$, with the convention that 
$\min \left( \emptyset \right)=N$.
Moreover
\begin{equation}
\label{eq:expm}
\bbE \, \bP^a_{N, \go}\left( \Ndown \ge m\right) \le \frac 1c \exp\left(-c m\right),
\end{equation}
both for $a = \rf$ and $a= \rc$, for every $N$ and every $m \in \N$. 
\medskip
\item
If the deviation inequality  holds
then 
for  $h>h_c(\gl )$ there exist two positive constants $c$ and $q$ such that
\begin{equation}
\bbE \, \bP^a_{N, \go}\left( \Ndown \ge m\right) \le  \exp\left(-c m\right),
\end{equation}
both for $a=\rf$ and $a= \rc$,
for every $N$ and every $m \ge q \log N$.
\end{enumerate} 
\end{theorem}

\medskip

We refer to Section \ref{sec:further}
for a thorough discussion on how these results relate
to what is expected to happen, with a particular attention
to scaling limits and almost sure results.
In the same section one finds also some further considerations
on the delocalized path behavior.

\smallskip
\begin{rem}\rm
The methods of proof of Theorem~\ref{th:main}
are applicable in more general contexts. We mention in particular
the case of disordered pinning or wetting.
Consider in particular the case of a model defined
like in \eqref{eq:Boltzmann}, but with $\sign(S_n)$  replaced
by $\ind_{\{0\}}\left(S_n \right)$. In spite of the formal resemblance,
this is a profoundly different model and, in order to deal with
interesting phenomena, one has to allow $h$ to take  negative values too.
In \cite{cf:AS} 
it is proven that the free energy of the model
exists and it is non negative and, exactly in analogy with $f(\gl,h)-\gl h$ in our setting,
one defines the localization and delocalization regimes  depending on whether the
free energy is positive or zero. Moreover for any $\gl>0$ the
transition takes place at a critical value $h_c$ of the parameter $h$
and $h_c \in [h_c^a, 0)$, where $h_c^a<0$ is the critical value
for the corresponding annealed model, a homopolymer model
that can be solved exactly. It is expected, but not proven, that
$h_c >h_c^a$ (see references in   \cite{cf:AS}).
    Theorem~\ref{th:main} holds for this disordered pinning model
 provided one changes the definition of $\cN$ to 
 $\cN:=\sum_{n=1}^N \ind_{\{0\}}\left(S_n \right)$.    
\end{rem}

\section{Proof of Theorem \ref{th:main}}
\label{sec:proofs}

\setcounter{equation}{0}

It is  convenient to consider
a modified partition function. To this purpose we observe that we may write
\begin{equation}
\frac{\dd \bP _{N, \go}^{a}}{\dd \bP}
\left(S \right)
\, = \, 
\frac 1{Z_{N, \go}^a}
\exp 
\left(-2 \gl \sum_{n=1}^{N} \left(\go_n +h\right) \Delta_n
\right) \ind_{\gO_N^a}\left( S\right),
\end{equation}
where $Z^a_{N, \go }= Z_{N, \go} (\gO_N ^a)$, with the notation 
\begin{equation}
 Z_{N,\go} \left( \tilde \gO \right) \, =\, \bE\left[
\exp\left(
-2 \gl \sum_{n=1}^{N} \left(\go_n +h\right) \Delta_n
\right); \, \tilde\gO
\right],
\end{equation}
for $ \tilde \gO\subset \gO$.
Likewise 
we introduce $\tf _{N, \go} \left(  \tilde \gO \right)
:= (1/N)\log Z_{N,\go} \left(  \tilde \gO \right)$ and
$\tf ^a_{N , \go}  := 	 \tf _{N , \go} ( \gO_N ^a)$.
Notice that $\bbP (\dd \go)$--a.s.
we have that $ Z_{N, \go}^a \asymp \tilde Z^a_{N,\go} \exp(-\gl h N )$, where
$\asymp$ denotes the Laplace asymptotic equivalence,
which means that the $\bbP(\dd \go)$--a.s. limit 
of $\tf ^a_{N , \go} $ equals $f(\gl, h)-\gl h=: \tf (\gl ,h)$.

\smallskip

\subsection{The concentration lemma}
For $m\in \{0,2, \ldots, N\}$ let us consider an event
$\gO_m \subset \gO$ such that $\bP \left( \gO_m \right)>0$
and such that $\Ndown  =m$ for every  $S\in \gO_m$.
If the distribution of $\omega$ satisfies the deviation inequality, we have

\medskip

\begin{lemma}
\label{th:concentration}
For every $N$, every $m \in \{0,2, \dots, N\}$ and every $u \ge 0$ we have
\begin{equation}
\bbP\left(
\tf_{N,\go} (\gO_m) -  \bbE\left[ \tf_{N,\go} (\gO_m) \right]\ge u
\right) \le C \exp\left(- \frac{u^2 N^2}{4C \gl^2 m}\right). 
\end{equation}
\end{lemma}

\medskip
\noindent
{\it Proof of Lemma \ref{th:concentration}}.
By the deviation inequality it suffices to show that for every $\omega,\omega'\in {\mathbb R}^N$ 
we have
\begin{equation}
\label{eq:Lip1}
\left\vert  \tf_{N, \go} \left(\gO_ m\right) - \tf_{N, \go^\prime} \left(\gO_ m\right)
\right \vert \, \le \, \frac{2\gl \sqrt{m} }{N} \left\Vert \go -\go^\prime \right\Vert,
\end{equation}
where $\Vert \,  \cdot \, \Vert$ is the Euclidean norm of $\cdot$.
In order to establish \eqref{eq:Lip1} we introduce $\go_t = t \go +(1-t)\go ^\prime$ 
and, taking the derivative with respect to $t$ and integrating back,  after the use
of the Cauchy--Schwarz inequality we obtain
\begin{equation}
\begin{split}
\left\vert  \tf_{N, \go} \left(\gO_m\right) - \tf_{N, \go^\prime} \left(\gO_m\right)
\right \vert & = \left\vert  \int_0^1 \sum_{n=1}^N \frac{(-2\gl)}N 
\bE_{N,  \go_t}\left[ \gD_n \big\vert \gO_m \right]
\left( \go_n - \go_n^\prime\right) \dd t \right\vert
\\
& \le \frac{2\gl}N \sqrt{ \sup_t \sum_{n=1}^N (\bE_{N,  \go_t}\left[ \gD_n \big \vert \gO_m\right]
)^2 
\sum_{n=1}^N \left( \go_n -\go_n^\prime \right)^2 }.
\end{split}
\end{equation}
Since $\sum_{n=1}^N (\bE_{N,  \go_t}\left[ \gD_n | \gO_m\right]
)^2 
\le m$,
 the proof is complete.
$\stackrel{\text{Lemma } \ref{th:concentration}}{\qed}$

\subsection{The delocalized  region}

\noindent
{\it Proof of Theorem \ref{th:main}, part (2).}
In this proof we set $\tf^a_{N, \go}(\gl, h):=(1/N) \log Z^a_{N, \go}$ and
\begin{equation}
\label{eq:213}
\tf^a_{N, \go}(\gl, h;m):=\frac 1N \log Z_{N, \go}\left(\gO_N^a\cap\left\{\Ndown =m\right\}\right).
\end{equation}
Since for $(\gl, h) \in \cD$ we have $\tf (\gl ,h)=0$ and since $\left\{ N\bbE\left[ \tf^\rc_{N, \go} (\gl, h)\right]\right\}_N$
is superadditive, so that $\lim_{N \to \infty} \bbE\left[\tf^\rc_{N, \go} (\gl, h)\right] =
\sup_N \bbE\left[\tf^\rc_{N, \go} (\gl, h)\right]$, we have that 
\begin{equation}
\label{eq:indeloc}
\bbE\left[ \tf^\rc_{N, \go}(\gl,h)\right] \le 0, 
\end{equation}
for every $N$. The superadditivity is a direct consequence of the
Markovian character of $S$, see \cite{cf:BdH} or \cite{cf:G} for the details.

Now let us fix $(\gl, h ) \in \stackrel{\circ}{\cD}$ and $\gep>0$ such that
$(\gl, h-\gep) \in \cD$. Observe that for every $\go$ 
\begin{equation}
\label{eq:equival}
\tf^\rc_{N,\go} (\gl,h; m)\ge -\gl \gep m/N \ \ \Longleftrightarrow \ \
\tf_{N,\go}^\rc (\gl,h-\gep ; m)\ge \gl \gep m/N,
\end{equation} 
but $\bbE\left[ \tf^\rc_{N, \go}(\gl,h-\gep; m)\right] \le \bbE\left[ \tf ^\rc_{N, \go}(\gl,h-\gep)\right] 
\le 0$, so that, by Lemma \ref{th:concentration}, we have
\begin{equation}
\label{eq:shifth}
\begin{split}
\bbP\left( \tf ^\rc_{N,\go}(\gl, h; m) \ge -\gl \gep m/N \right) 
&=
\bbP\left(\tf ^\rc_{N,\go} (\gl,h-\gep ; m)\ge \gl \gep m/N\right) 
\\ 
&\le C\exp\left(-  \gep^2 m/4C\right).
\end{split}
\end{equation}
From this we directly obtain that if we set
 $E_{\overline{m}}= \{$there exists $m\ge \overline{m}$ such that 
 $\tf ^\rc_{N,\go}(\gl, h; m) \ge -\gl \gep m/N\}$ 
then
\begin{equation}
\label{eq:pe}
\bbP\left(E_{\overline{m}} \right) \le
c_1 \exp\left(- c_2  \overline{m}\right).
\end{equation}
We can now evaluate the tail of $\Ndown$.
For  $\go \in E_{\overline{m}}^\complement$, with $\gO_m = \{\cN =m, \, S_N=0\}$, we have
\begin{equation}
\label{eq:mw1}
\bP_{N,\go }^\rc
\left( \Ndown \ge \overline{m} \right) \, = \, 
\frac{\sum_{m \ge \overline{m}}  Z_{N, \go}\left(\gO_m\right)}{
Z^\rc_{N, \go}} \, \le \,  c_3 N^{3/2} \sum_{m \ge \overline{m}} 
\exp\left(-\gl \gep m\right) 
\, \le \,  c_4 N^{3/2} \exp\left(-\gl \gep \overline{m}\right)  ,
\end{equation}
where we have used
that $ Z^\rc_{N, \go}\ge  \bP (S_n>0, \, n=1, \ldots , N-1, \, S_N=0)\ge 1/(c_3 N^{3/2})$
\cite[Ch. 3]{cf:Feller1}.
The estimates
\eqref{eq:pe} and
\eqref{eq:mw1} readily imply
\begin{equation}
\bbE \, \bP^\rc_{N, \go}\left( \Ndown \ge \overline{m} \right) \le c_5
N^{3/2} \exp\left(-c_6 \overline{m}\right).
\end{equation}
The choice of $\overline{m}\ge q \log N$, for $q$ sufficiently large completes 
the proof for the case of $ \bP^\rc_{N, \go}$.

\smallskip

For the free endpoint case $\bP^\rf_{N, \go}$ 
one recalls that in  \cite{cf:BdH} (or in \cite{cf:G})
it is proven that there exists a positive constant $c$ such that 
\begin{equation}
\label{eq:fromBdHandG}
Z^\rf_{N, \go } \le  c N  Z^\rc_{N, \go },
\end{equation}
for every $\go$ and every $N$.
Therefore, by \eqref{eq:indeloc}, we have
\begin{equation}
\bbE\left[ \tf^\rf_{N, \go}(\gl,h)\right] \, \le\, \frac 1N \log (cN), 
\end{equation}
and therefore formulas \eqref{eq:equival}, \eqref{eq:shifth} and  \eqref{eq:pe}
hold if we replace the quantity $\tf^\rc_{N, \go}(\gl,h; m)$
with $\tf^\rf_{N, \go}(\gl,h; m)- ( \log cN)/N$.
It suffices therefore to observe that
$ \inf_{N, \go} N^{1/2}Z^\rf_{N, \go}\ge  \inf _N N^{1/2}
\bP (S_n>0, \, n=1, \ldots , N)>0$
to conclude that
\eqref{eq:mw1} holds unchanged if $a=\rc$ is replaced by $a=\rf$
and $\gO_m= \left\{\cN =m \right\}$,  
apart for the explicit values of the multiplicative constants (which we have not tracked anyway). 
The proof is therefore complete.
$\stackrel{\text{Theorem }\ref{th:main}(2)}{\qed}$

\subsection{The strongly delocalized region}
\noindent
{\it Proof of Theorem \ref{th:main}, part (1).}
We start by proving \eqref{eq:expm}. 
We compute 
by means of the 
Fubini--Tonelli theorem:
\begin{equation}
\label{eq:QtoA}
\begin{split}
\bbE\left[  Z_{N, \go} \left( \gO_m \right)\right] \, & = \, 
  \bE\bbE\left[
\exp\left( -2 \gl \sum_{n=1}^N \left(\go_n+h \right) \Delta_n  \right)
; \, \gO_m \right]
\\
& =\,
 \bE\left[
\exp\left( \sum_{n=1}^N \left( \log M( 2\gl \gD_n) -2\gl h \gD_n \right)\right); \, \gO_m \right]
\\
& =\,
 \bE \left[\exp  \left(
-2\gl \left( h - \log M( 2\gl )/ 2\gl \right) 
\sum_{n=1}^N \gD_n \right) \bigg \vert \gO_m 
\right] \, 
 \bP \left( \gO_m
\right),
\end{split}
\end{equation}
where $\gO_m$ is like in Lemma \ref{th:concentration}.
Since $\Ndown = \sum_{n=1}^N \gD_n =m$ on $\gO_m$ we have
\begin{equation}
\label{eq:QtoA1}
\bbE\left[  Z_{N, \go}\left( \gO_m \right) \right] \, 
= \, 
\bP \left( \gO_m
\right) \exp\left( -\gb m\right),
\end{equation}
with $\gb :=  2\gl \left( h - \log M( 2\gl )/ 2\gl \right)$, so $\gb>0$ in the strongly
delocalized regime.

\smallskip
We can now estimate the tail behavior of  $\Ndown $,
averaged over the disorder $\go$. 
We first consider the free case: set
$\gO_m = \{  \Ndown =m\}$ and  $P_N(m):=\bP \left( \gO_m \right)$.
We have
\begin{equation}
\label{eq:mw2}
\bbE\, \bP_{N,\go }^\rf
\left( \Ndown \ge \overline{m} \right) \, = \, \bbE \left[
\frac{\sum_{m \ge \overline{m}}  Z_{N, \go}\left(\gO_m\right)}{
Z^\rf_{N, \go}}\right] \, \le \sum_{m \ge \overline{m}} 
\exp\left(-\gb m\right)
\frac{P_N(m)}{P_N(0)} ,
\end{equation}
where we have used once again that 
$Z_{N, \go}^\rf \ge \bP \left( \cN =0\right) =P_N(0)$ and we have applied \eqref{eq:QtoA1}.
The proof of \eqref{eq:expm} in the free endpoint case is completed once we observe
that $ P_N(m) \le P_N(0)$, a fact that can be easily extracted from the exact
expression of $P_N(m)$  \cite[Ch. 3]{cf:Feller1}.

In the constrained endpoint case one takes $\gO_m= \left\{ \Ndown =m, \, S_N=0\right\}$ 
and the steps are then identical. Notice however that in this case 
$P_N(m)=P_N (0)$ every $m$  \cite[Ch. 3]{cf:Feller1}.  

\bigskip 

We turn now to the proof of \eqref{eq:withoutconc1} and  \eqref{eq:withoutconc2}.
Like in \eqref{eq:QtoA}, by explicit computation we have
\begin{equation}  
\bbE \left[  Z^a_{N, \go}\right] \, =\,
\bE\left[ \exp\left(-\beta\Ndown \right); \, \gO_N^a\right]. 
\end{equation}

Let us observe preliminarily that we have  \cite[Ch. 3]{cf:Feller1}:
\begin{equation}
\label{eq:fromFeller1}
\bP\left( \Ndown= k\right) \le \frac c{\sqrt N} \  \  \  \text{ and } \ \ \ 
\bP\left( \Ndown = k,  S_{N}=0\right)\le \frac c{ N^{3/2}},
\end{equation}
and from this one easily finds a constant $C>0$ such that
\begin{equation}
\label{eq:UBnum}
\bbE \left[  Z^\rf_{N, \go}\right] \le\frac{C}{N^{1/2}}
 \  \  \  \text{ and } \ \ \ 
\bbE \left[   Z^\rc_{N, \go}\right] \le\frac{C}{N^{3/2}}.
\end{equation}
On the other hand we know (and used several times by now) that there exists $c>0$ such that
\begin{equation}
\label{eq:LBden}
Z^\rf _{N, \go} \, \ge\,  c/ N^{1/2} \ \ \ \text{ and } \ \ \ 
Z^\rc _{N, \go} \, \ge\,  c/ N^{3/2}
\end{equation}
for every $\go$.

\smallskip
We focus now on the proof of \eqref{eq:withoutconc1}. 
Let us call $F_\ell$ the event of the random walk trajectories 
for which there exists $n \in \{\ell, \ldots ,N\}$  such that $S_{\ell}=0$.
By conditioning on the last hitting time of zero before time $N$ we obtain
\begin{multline}
\bP^\rf_{N, \go} \left(F_\ell\right)
\, = \, \\
 \frac 1{Z_{N,\go}^\rf} \sum_{l =\ell}^N  Z^\rc_{l, \go}\bP\left(
S_n>0 \text{ for } n=1, 2, \ldots, (N-l)
\right) 
\left(1+ \exp\left(-2\gl \sum_{n=l+1}^N (\go_n +h)\right) \right). 
\end{multline}
Since the denominator can be bounded below uniformly in $\go$, cf.
\eqref{eq:LBden}, by integrating with respect to $\go$ we obtain
\begin{multline}
 \bbE\left[
 \bP^\rf _{N, \go} \left(F_\ell \right)
 \right]\, \le \,
 c  N ^{1/2} \sum_{l = \ell}^N  \bbE\left[  Z^\rc_{l, \go} \right]
 \bP\left(
S_n>0,\,  n=1, 2, \ldots, N-l)\left( 1+ \exp(-\gb (N-l)\right)
\right)
\\
\le  c_1 N^{1/2}\sum_{l = \ell}^N \frac1 {(l+1)^{3/2}} \frac1 {(N+1-l)^{1/2}}
\, \le\, 
\frac{c_2}{(\ell+1)^{1/2}}.
\end{multline}
In order to complete the proof of \eqref{eq:withoutconc1}
we need to exclude that the last excursion of the polymer is in the lower half plane.
However one directly verifies that:
\begin{equation}
\label{eq:downexcurs}
\bbE \,\bP^\rf _{N \go}\left( F_\ell ^\complement, \, S_n <0 \text{ for }n \ge \ell\right) \le
\exp( -\gb (N-\ell)),
\end{equation}
and this suffices to conclude the proof of  \eqref{eq:withoutconc1}.

\smallskip

The proof of 
\eqref{eq:withoutconc2} is conceptually very close to the proof 
of \eqref{eq:withoutconc1}.
We introduce the event $F_{\ell_1, \ell_2}$
of the polymer trajectories hitting $0$ 
in the set $\{\ell_1,\ldots, N/2 \}$
and in $\{N/2,\ldots, N -\ell_2 \}$. We may write 
\begin{equation}
\bP_{N, \go }^\rc\left( F_{\ell_1, \ell_2}\right)
\, =\, \frac 1{Z_{N, \go}^\rc}
\sum_{j_1=\ell_1}^{N/2} \sum_{j_2= \ell _2}^{N/2}
Z_{j_1,  \go}^\rc Z^\pm_{N, \go}(j_1,j_2)
Z_{j_2,  \tau_{N-j_2}\go}^\rc,
\end{equation}
where  $\tau_k$ is the $k$--shift, i.e., $(\tau_k \go)_n= \go_{n+k}$,
and  
\begin{multline}
Z^\pm_{N, \go}(j_1,j_2) \, :=\, 
\bP \left( S_n>0, \, n=1,2, \ldots, N-j_1-j_2-1, \, S_{ N-j_1-j_2}=0\right)\\
\times
\left( 1+ \exp\left( -2\gl \sum_{n=j_1+1}^{N-j_2} (\go_n +h )\right) \right).
\end{multline}
Once again we estimate the denominator uniformly with respect to
$\go$, cf. \eqref{eq:LBden}, and then take expectation. By applying 
\eqref{eq:fromFeller1} and \eqref{eq:UBnum} we obtain
\begin{equation}
\bbE\, \bP_{N, \go }^\rc\left( F_{\ell_1, \ell_2}\right) \, \le \, 
c_1 \sum_{j_1= \ell_1}^{N/2}
\sum_{j_2= \ell_2}^{N/2}
\frac1{(j_1+1)^{3/2}}
\frac1{(j_2+1)^{3/2}}
\left( \frac N{N-j_1-j_2 +1}\right)^{3/2} .
\end{equation}
Since the right--hand side can be bounded above by $c_2/\sqrt{\ell _1 \ell_2+1}$
the proof is easily completed. 
$\stackrel{\text{Theorem }\ref{th:main}(3)}{\qed}$
\bigskip

\section{Universality of the slope at the origin}
\label{sec:exapp}

{\em Proof of Theorem \ref{univers}.}
Let us first of all prove the theorem 
 when the random variable $\omega_1$ is bounded.

Let ${\mathbb P}^{(1)}$ be the law of IID centered Bernoulli random variables $\omega_n=\pm1$.
Also, consider a law ${\mathbb P}^{(2)}$ for the
IID symmetric bounded random variables $\{\omega_n\}_n$, 
and recall that by convention 
${\mathbb E}^{(2)} [\omega_1^2]=1$.

Let $m_c^{(1)}$ be the slope at the origin
of
the critical curve of the  copolymer model with Bernoulli disorder, whose existence was proven in
\cite{cf:BdH}. By definition of $m_c^{(1)}$, for any $v>m_c^{(1)}$ and 
$\lambda$ sufficiently small one has 
\begin{eqnarray}
  \tf^{(1)}(\lambda,v\lambda):= \sup_{N}
{\mathbb E}^{(1)} \tf^\rc_{N,\omega}(\lambda,v\lambda)=0.
\end{eqnarray}
This implies that, for any $\varepsilon>0$, $N\in 2{\mathbb N}$ and
$m\in\{0,2,\cdots,N\}$,
\begin{eqnarray}
\label{eq:epsilon}
{\mathbb E}^{(1)} \tf^\rc_{N,\omega}(\lambda,(v+
\varepsilon)\lambda;m)\le -\frac{2\varepsilon \lambda^2 m}N,
\end{eqnarray}
where we use the same notation as in equation (\ref{eq:213}).
On the other hand, one has the following lemma, proven below.

\medskip

\begin{lemma}
\label{lemmicchio}
If the laws ${\mathbb P}^{(\ell)}$, $\ell=1,2$ correspond to IID centered bounded random variables,
there exists a constant $c>0$ such that for any $h\ge0$, $0\le \lambda\le1$, $N\in 2{\mathbb N}$ and
$m\in\{0,2,\cdots,N\}$,
  \begin{eqnarray}
\label{eq:interp}
    \left|{\mathbb E}^{(1)} \tf^{a}_{N,\omega}(\lambda,h;m)
-{\mathbb E}^{(2)} \tf^{a}_{N,\omega}(\lambda,h;m)\right|\le c\frac{m\lambda^3}N,
  \end{eqnarray}
  and
  \begin{eqnarray}
  \label{eq:geninterp}
    \left|{\mathbb E}^{(1)} \tf^{a}_{N,\omega}(\lambda,h)
-{\mathbb E}^{(2)} \tf^{a}_{N,\omega}(\lambda,h)\right|\le c\lambda^3.
  \end{eqnarray}
\end{lemma}
\medskip

Thanks to equations (\ref{eq:epsilon}) and (\ref{eq:interp}), one has 
\begin{eqnarray}
\label{eq:unifbound}
{\mathbb E}^{(2)}\tf^\rc_{N,\omega}(\lambda,(v+
\varepsilon)\lambda;m)\le -\frac{2\varepsilon \lambda^2
  m}N+c\frac{m\lambda^3}N
\le -\frac{\varepsilon \lambda^2m}N
\end{eqnarray}
provided that $\lambda\le\min(1,\varepsilon/c)$.
Using the deviation inequality (\ref{eq:Lip}), which is applicable since the random variables are bounded,
it is then possible to deduce that
\begin{eqnarray}
\label{eq:spezz}
  \tf^{(2)}(\lambda,(v+\varepsilon)\lambda):= \lim_{N\to\infty}
{\mathbb E}^{(2)}\tf^{\rc}_{N,\omega}(\lambda,(v+\varepsilon)\lambda)=0.
\end{eqnarray}
This point is discussed in greater detail at the end of this section, in a more general context where
the random variables $\omega_n^{(2)}$ are not necessarily bounded.
Therefore, one has
$\limsup_{\lambda\searrow 0}{h_c^{(2)}(\lambda)}/\lambda\le v+\varepsilon$
and, thanks to the arbitrariness of 
$\varepsilon>0$ and of $v>m_c^{(1)}$,
\begin{equation}
\limsup_{\lambda\searrow 0}\frac{h_c^{(2)}(\lambda)}\lambda\le m^{(1)}_c.
\end{equation}
To obtain the opposite bound, observe that from Theorem 6 of \cite{cf:BdH} follows that
for any $v<m^{(1)}_c$, there exists $c(v)>0$ such that
\begin{eqnarray}
\tf^{(1)}(\lambda,v\lambda)\ge c(v)\lambda^2
\end{eqnarray}
for $\lambda$ sufficiently small. On the other hand,
thanks to \eqref{eq:geninterp}, for $\lambda$ sufficiently small one has
\begin{eqnarray}
\tf^{(2)}(\lambda,v\lambda)\ge \frac{c(v)}2\lambda^2
\end{eqnarray}
which implies 
\begin{equation}
\liminf_{\lambda\searrow 0}\frac{h_c^{(2)}(\lambda)}\lambda\ge m^{(1)}_c
\end{equation}
and the statement of the theorem in the bounded case.
$\stackrel{\text{Theorem }\ref{univers}, \ \go_1 \text{ bounded}}{\qed}$

{\em Proof of Lemma \ref{lemmicchio}} This is based on an interpolation 
argument, of the type of the one showing that the free energy
of the Sherrington-Kirkpatrick spin glass model does not 
depend on the distribution of the couplings (see \cite{cf:talacergy}, \cite{cf:limite} and the 
more recent \cite{cf:ch}).

For definiteness, we give the proof of \eqref{eq:interp} in the pinned case $a=\rc$.
For $0\le t\le1$, consider the auxiliary free energy
\begin{eqnarray}
  F_N(t)=\frac1N {\mathbb E}^{(1,2)} \log \bE\left[ \exp\left(-2\lambda
\sum_{n=1}^N(\sqrt t \omega^{(1)}_n+\sqrt{1-t} \omega^{(2)}_n+h)
\Delta_n\right);{\mathcal N}=m,S_N=0\right],
\end{eqnarray}
where $\omega^{(1)},\omega^{(2)}$ are independent and distributed according to 
the laws ${\mathbb P}^{(1)},{\mathbb P}^{(2)}$ respectively.
Then, one has immediately
\begin{eqnarray}
  && F_N(1)=\bbE^{(1)}\tf^{\rc}_{N,\omega}(\lambda,h;m)\\
&& F_N(0)=\bbE^{(2)}\tf^{\rc}_{N,\omega}(\lambda,h;m).
\end{eqnarray}
Therefore, one has to estimate the $t$-derivative of the free energy, which is easily computed: 
\begin{eqnarray}
\label{eq:derivata}
\frac {\dd F_N(t)}{\dd t}=-\frac\lambda N{\mathbb E}^{(1,2)}\sum_{n=1}^N \bE_{N,\omega_t}(\Delta_n|{\mathcal N}=m,S_N=0)
\left(\frac1{\sqrt t}\omega^{(1)}_n-\frac1{\sqrt{1-t}}\omega^{(2)}_n\right).
\end{eqnarray}
This expression can be manipulated by means of the identity
\begin{eqnarray}
\label{miparegiustamacontrolla}
\bbE \,\eta G(\eta)= \bbE\, G'(\eta)+\bbE\left((\eta^2-1)\int_0^\eta G''(u)\dd u\right)
-\frac14\bbE\, |\eta|\int_{-|\eta|}^{+|\eta|}(\eta^2-u^2)G'''(u)\dd u,
\end{eqnarray}
which holds for any symmetric random variable $\eta$  
with $\bbE\,[ \eta^2]=1$ and for sufficiently regular functions $G$. 
In our case, the idea is that every derivative with respect to $\omega_n^{(\ell )}$ carries a (small) factor
$\lambda$, so that the first term in the r.h.s. of (\ref{miparegiustamacontrolla}) is the 
dominant one. Applying this identity to (\ref{eq:derivata}), one finds that the {\sl dominant terms}
cancel exactly, and one is left with terms involving derivatives of $\bE_{N,\omega_t}(\Delta_n|{\mathcal N}=m,S_N=0)$
of order higher than one. Indeed, denoting 
$$
0\le X_n^{(\ell)}(u):= \bE_{N,\omega_t}(\Delta_n|{\mathcal N}=m,S_N=0)\big|_{\omega^{(\ell)}_n=u},
$$
and noting that for $k\ge1$
$$
0\le (X_n^{(\ell)}(u))^k\le X_n^{(\ell)}(u),
$$
one has
\begin{eqnarray}
\label{linea1}
\left|\frac {\dd F_N(t)}{\dd t}\right|&\le&
\sum_{\ell=1}^2\frac{12\lambda^3}N\sum_{n=1}^N\bbE^{(1,2)} ((\omega^{(\ell)}_n)^2+1)\left|
\int_0^{\omega^{(\ell)}_n}
X_n^{(\ell)}(u)\,\dd u\right|\\
&&\label{linea2}
+\sum_{\ell=1}^2
\frac{26\lambda^4}N\sum_{n=1}^N\bbE^{(1,2)} |\omega^{(\ell)}_n|^3\int_{-|\omega^{(\ell)}_n|}^{+|\omega^{(\ell)}_n|}
X_n^{(\ell)}(u)\,\dd u.
\end{eqnarray}
Below, we consider only the terms with $\ell=1$, the other case requiring only minimal modifications.
Let us first consider the term in (\ref{linea2}). Observe that
\begin{eqnarray}
-2\lambda X^{(1)}_n(u)\le -2\lambda\sqrt t\left(X^{(1)}_n(u)- (X^{(1)}_n(u))^2 \right)=\frac d{du}X^{(1)}_n(u)\le0
\end{eqnarray}
so that for any $u,u'$
\begin{equation}
X^{(1)}_n(u)\le X^{(1)}_n(u')e^{2\lambda|u-u'|}.
\end{equation}
Therefore, the term in (\ref{linea2}) can be bounded above by
\begin{equation}
\frac{26\lambda^4}N\sum_{n=1}^N \bbE^{(1,2)} |\omega^{(1)}_n|^3\int_{-|\omega^{(1)}_n|}^{+|\omega^{(1)}_n|} \,
\bE_{N,\omega_t}(\Delta_n|{\mathcal N}=m,S_N=0)e^{2\lambda|u-\omega^{(1)}_n|}\dd u\, \le\, c\frac{m\lambda^4}N
\end{equation}
where we made use of the boundedness of $\omega_n^{(1)}$ and of the fact that
$\bE_{N,\omega_t}(\Delta_n|{\mathcal N}=m,S_N=0)$ does not depend on $u$.
An analogous bound can be obtained for the term in (\ref{linea1}). Indeed, it is bounded above by
\begin{multline}
c\frac{\lambda^3}N\sum_{n=1}^N\bbE^{(1,2)} 
\left|\int_0^{\omega^{(1)}_n}\bE_{N,\omega_t}(\Delta_n|{\mathcal N}=m,S_N=0)
e^{2\lambda|u-\omega^{(1)}_n|}\,\dd u\right|\\
\le c'\frac{\lambda^3}N\sum_{n=1}^N\bbE^{(1,2)}\bE_{N,\omega_t}(\Delta_n|{\mathcal N}=m,S_N=0)=
c^\prime \frac{m\lambda^3}N.
\end{multline}
The proof of \eqref{eq:geninterp} is much simpler. In this case one removes the constraint on $\mathcal N$ in the 
definition of $F_N(t)$ and then it is immediate to realize that \eqref{linea1} and \eqref{linea2} are of 
order $O(\lambda^3)$ and $O(\lambda^4)$, respectively.

$\stackrel{\text{Lemma } \ref{lemmicchio}}{\qed}$

\medskip

\noindent
{\it Proof of Theorem \ref{univers}, $\go_1 \sim \cN (0,1)$.}
It remains to show that the proof covers also the case when ${\mathbb P}^{(2)}$ is the law of IID centered 
Gaussian variables. One easily verifies 
that Lemma \ref{lemmicchio} still holds if one of the two laws is replaced by 
the Gaussian one. To this purpose, it is
sufficient to observe that, if $\eta$ is a ${\mathcal N}(0,1)$ random variable, identity 
(\ref{miparegiustamacontrolla}) can be replaced by the integration by parts formula
$$
\bbE\,\eta G(\eta)= \bbE \,G'(\eta).
$$
Therefore, one still obtains the uniform bound (\ref{eq:unifbound}).
In order to deduce (\ref{eq:spezz}) from (\ref{eq:unifbound}), one proceeds as follows. 
For any $\overline m\ge0$, one can decompose the partition function and write with obvious notation
\begin{equation}
Z^\rc_{N,\omega}(\lambda,h)=Z^\rc_{N,\omega}(\lambda,h;m\le\overline m)+Z^\rc_{N,\omega}(\lambda,h;m>
\overline m).
\end{equation}
From now until the end of the proof we set $h=(v+\varepsilon)\lambda$.
Then, using the inequality
\begin{equation}
\log(a+b)\le \log 2+\log a+\log b,
\end{equation}
which holds whenever $a,b\ge1$, and the fact that $Z^\rc_{N,\omega}(\lambda,h;m\le\overline 
m)\ge c_1 N^{-3/2}$ for
some constant $c_1$ independent of $\omega$ and $N$, one has
\begin{eqnarray}
\label{eq:pdCS}
\bbE^{(2)} \tf^\rc_{N,\omega}(\lambda,h)&\le& \frac1N\bbE^{(2)}\log Z^\rc_{N,\omega}(\lambda,h;m\le
\overline m)
\\\nonumber
&&+\frac1N\bbE^{(2)}\log \max\left(1,Z^\rc_{N,\omega}(\lambda,h;m>\overline m)N^{3/2}/c_1\right)
+\frac{\log 2}N.
\end{eqnarray}
The first term in (\ref{eq:pdCS}) can be bounded above via Jensen's inequality by
$ {c_2\overline m}/N$, where $c_2$ is a constant independent of $\overline m$ and $N$.
As for the second term, define the event $E_{\overline m}$ as
\begin{equation}
E_{\overline m}=\left\{\mbox{ there exists\;\;} m\ge\overline m \mbox{ such that\;\;} 
\tf^\rc_{N,\omega}(\lambda,h;m)\ge -\frac{\varepsilon m \lambda^2}{2N}\right\}
\end{equation}
whose probability, thanks to the deviation inequality (\ref{eq:Lip}) and to \eqref{eq:unifbound}, satisfies
$$
\bbP^{(2)} \left (E_{\overline m}\right)\le \frac1{c_3} e^{-c_3\overline m}.
$$
The second term in  (\ref{eq:pdCS}) can be therefore bounded above by
\begin{eqnarray}
c_4\frac{\log N}N+\frac1N\sqrt{\bbP^{(2)} \left (E_{\overline m}\right)\bbE^{(2)}\left(
\log\max\left(1,Z^\rc_{N,\omega}(\lambda,h;m>\overline m)N^{3/2}/c_1\right)\right)^2
},
\end{eqnarray}
where the first term comes from the average restricted to the event 
$E_{\overline m}^\complement$ and in the second 
we applied Cauchy-Schwarz inequality.
Observing that 
$$
Z^\rc_{N,\omega}(\lambda,h;m>\overline m)\le \exp\left({2\lambda\sum_{n=1}^N(h+|\omega_n|)}\right)
$$ 
and 
putting everything together, one obtains finally
\begin{eqnarray}
\bbE^{(2)} \tf^\rc_{N,\omega}(\lambda,h)&\le&c_5\left(\frac {\overline m}N+
\frac {\log N}N+e^{-c_3 \overline m/2}\right),
\end{eqnarray}
from which (\ref{eq:spezz}) follows choosing for instance $\overline m=\sqrt N$.

From this point on, the proof proceeds exactly like in the bounded case.
$\stackrel{\text{Theorem }\ref{univers}, \ \go_1\sim \mathcal N(0,1)}{\qed}$

\section{Further results and considerations}
\label{sec:further}

\setcounter{equation}{0}

\subsection{What does one expect on delocalized paths}
In the previous section, we have given delocalization results that hold 
in average with respect to the $\bbP$-probability. On the other hand, one would like to prove 
$\bbP(\dd \omega)$--almost sure results. 
In this respect, based on what is known  on non-disordered
models, see e.g. \cite{cf:DGZ} and \cite{cf:IY}, 
it is tempting to conjecture the following scenario:
for $h>h_c(\gl)$
\smallskip

\newcounter{Lcount}
\begin{list}{C.\arabic{Lcount}}
      {\usecounter{Lcount}
      \setlength{\rightmargin}{\leftmargin}}
    \item \label{C1} there are only a finite number of visits to the unfavorable solvent, that is 
    $\bbP( \dd \go )$--a.s.
    \begin{equation}
    \label{eq:C1}
    \lim_{\ell \to \infty} \limsup_{N\to \infty}
    \bP^\rf _{N, \go}
    \left( \max \cA  >\ell\right) =0.
    \end{equation} 
    \item \label{C2} there is a diffusive scaling limit to a {\sl Brownian meander}. In other terms
    if  we set $B_t ^{(N)}:= S_{Nt}/\sqrt{N}$ for $N \in \{0,1, \ldots, N\}$ and we extend the definition of $B_\cdot^{(N)}$ to a function in $C^0([0,1]; \R)$ by linear interpolation for $S_{N\cdot}$,
    $\bbP (\dd \go)$--a.s. 
we have that the law of  $B_\cdot^{(N)}$, with $S$
    distributed according to   $\bP_{N, \go}^\rf  $, converges 
    weakly as $N \to \infty$ to the law of the Brownian meander, that is the law of a standard Brownian
    process conditioned not to enter the lower half plane. 
The standard reference
    for the Brownian meander is \cite{cf:RY}.
  \end{list}
  
\smallskip

With the same level of confidence one might formulate the analogous 
statements for the constrained case: in particular, in C.2 the expected scaling limit
would the Brownian bridge conditioned to stay positive, a process that normally goes
under the name of {\sl Bessel(3) bridge}, see \cite{cf:RY}.

\smallskip 

As we will see, the  scenario outlined above 
cannot hold if taken literally, though we expect the 
qualitative picture to be correct.
To start with, observe that Theorem \ref{th:main} gives partial support to the conjectures, 
at least for $h\ge \overh(\lambda)$.  Indeed, for example if we choose a sequence $\{\ell_N \}_N$ with
$\lim_N \ell_N=\infty$, then  
by 
\eqref{eq:withoutconc1} we have
that 
\begin{equation}
\lim_N \bP ^\rf_{N, \go} \left( \max \cA > \ell_N\right)\, =\, 0,
\end{equation}
in $\bbP$--probability, or $\bbP( \dd \go )$--a.s. by subsequences.
This of course falls a bit short of proving C.1, even in the strongly delocalized region.
Just about the same is true for C.2. Let us set $\zeta_N := \max \cA$.
By the result we just stated, for $h> \overline{h}(\gl)$
there exists  a sequence $\left\{ N_j \right\}_j$ such that
$\bbP (\dd \go)$--a.s.  the random variable $\zeta_{N_j}/N_j$
   tends to
zero as $j  $ tends to infinity, in $\bP_{N_j, \go }$--probability.
Since it is not difficult to see that, conditionally
to $\zeta_N= k$, the law of $\left\{ S_{\zeta_N 
+n}\right\}_{n=0,1,\ldots}$
coincides with the law of a simple random walk constrained
not to enter the lower half--plane up to time $N-k$,
we are in the framework already considered for example in
\cite{cf:IY} or \cite{cf:DGZ}. Therefore, by proceeding like in  
\cite{cf:DGZ},
   one can show that for
a $\bbP$--typical $\go$ the sequence of random functions
  $\left\{ B ^{(N_j)}_\cdot\right\}_j$,
with $S$
     distributed according to   $\bP_{N_j, \go}^\rf  $, converges
     weakly as $j \to \infty$ to the law of the Brownian meander.

\smallskip

On the other hand, Theorem \ref{th:main} does not say much in the direction
 of C.1 and C.2 for $h \in (h_c(\gl), \overh (\gl)]$. This is due to the fact
 that, in spite of knowing that there are few visits to the unfavorable solvent,
 we do not know that they are close to the origin (or to $N$, in the constrained case).

\subsection{On the size of $Z_{N,\go}^a$}
Some further insight on the behavior of paths in the delocalized phase may 
be obtained by looking
at the size of $Z_{N , \go}^a$. 

Observe that, 
by \eqref{eq:UBnum}, under the basic assumptions on the disorder distribution
and in the strongly delocalized regime $h\ge \overh(\lambda)$,  $Z_{N,\omega}^a$ is of the
order of $N^{-1/2}$ for $a =\rf$, and $N^{-3/2}$ for $a =\rc$,
in the evident $\bbP$--probability sense. Recalling $\eqref{eq:LBden}$,
the result is somewhat sharp. 
We lack however an almost sure result
going beyond the fact that $Z_{N,\go} ^a$ tends to $0$
$\bbP(\dd \go) $--a.s. in the strongly delocalized region
(this is an immediate consequence of the convergence in probability 
along with the fact that $\{ Z_{N, \cdot}^a\}_N$ is a positive supermartingale
for $h \ge \overh (\gl)$
with respect to the natural filtration of $\go$ and therefore it
converges $\bbP (\dd \go)$--a.s.).

\smallskip

On the other hand, the result that we are going to present now says that something qualitatively
different happens for $h_c(\lambda)\le h\le \overh(\lambda)$. As we will discuss
at the end of the section, this phenomenon reflects on the behavior of the 
paths.

 \medskip
 
 \begin{proposition}
 \label{th:LBonZ}
 Under the basic assumptions on $\go$
 one can construct  a sequence $\left\{\tau_N \right\}_N$ of stopping times, with respect to the
 natural filtration of the sequence $\go$,  
 with the property that $ \log \tau_N (\go)/ \log N \stackrel{N \to \infty}{\longrightarrow} 1$
 $\bbP (\dd \go ) $--a.s.  and we can find a  number 
$\delta= \delta(\gl, h)$, explicitly given below, such that $ \delta >0$ if $h< \overh(\gl)$,
and  that
\begin{equation}
\label{eq:delta}
\lim_{N\to \infty} N^{1/2-\delta^\prime }Z_{\tau_N(\go), \go}^\rf =+\infty,
\ \ \ \ \bbP (\dd \go )-\text{a.s.}.
\end{equation}
for every $\delta^\prime < \delta$.
 \end{proposition}
 
 \medskip

\noindent
{\it Proof.}
Set $\tilde \go_n := \go_n +h$, choose a real number $q<h$
 and define
\begin{equation}
\tau_N \, := \,
\inf 
\left\{  n \in 2 \N : \frac{\sum_{j=k+1}^n \tilde \go _j}{n-k} \le  q \text{ for some } k \in\{0, 2, \ldots, n- r_N \}
 \right\}
 \end{equation}
with
$r_N$ the largest even integer smaller than $(\log N )/ \Sigma_h(q)$, where $\Sigma_h(\cdot)$
is the Cramer Large Deviation functional of  $\tilde \go$: 
\begin{equation}
\lim_{\ell \to \infty} \frac 1\ell \log \bbP\left( \sum _{j=1}^\ell \tilde \go_j \le  q\ell\right)\, =\, 
-\Sigma_h(q).
 \end{equation}
 We note that $\tau_N$ is the first moment at which an {\sl atypical } stretch of length
 at least $r_N$ appears along the
 sequence $\go$. 
 By Theorem 3.2.1 in \cite[\S 3.2]{cf:DZ}
 we have that
 $ \log \tau_N / \log N $ tends to $1$ 
 $\bbP (\dd \go ) $--a.s..
The same theorem tells us that, if
 \begin{equation}
 \label{eq:Rn}
 R_n \, :=\, \max \left\{
 \ell -k : \,  k \text{ and } \ell \text{ even }, \,   0 \le k < \ell \le
n , \,
 \frac{\sum_{j=k+1}^\ell \tilde \go _j}{\ell -k }\le q
 \right\},
 \end{equation}
 then $R_n/\log n \stackrel{n\to \infty}{\longrightarrow}  1/ \Sigma_h(q)$
$\bbP( \dd \go )$--a.s.
 and therefore
 \begin{equation}
 \lim_{N \to \infty} \frac{R_{\tau_N}}{\log N}\, =\, \frac1{\Sigma_h(q)},
\ \ \
 \bbP (\dd \go )-\text{a.s.}.
 \end{equation}
 Notice that the longest atypical stretch, in the sense of \eqref{eq:Rn},
 for $n=\tau _N$ ranges from $\tau_N - R_{\tau_N}$ to $\tau_N$,
 so $ \sum_{j=\tau_N - R_{\tau_N}+1}^{\tau_N} \tilde \go_j \le q
R_{\tau_N}$.

 \smallskip
 
 Choose now any $\gep>0$ and  a typical $\go$. We have
 \begin{equation}
 \begin{split}
 Z_{\tau_N ,\go}^\rf \, 
 & \ge \, \bP\left( S_n>0, \, n=1,2, \ldots, \tau_N  -R_{\tau_N} -1, \, 
 S_{ \tau_N -R_{\tau_N}} =0\right)
 \\
& \phantom{movemov} \times 
 \exp\left(-2\gl q R_{\tau_N}\right) \bP\left( S_n<0, \, n=1,2, \ldots,
 R_{\tau_N} \right)
 \\
 & \ge \, c \frac1{N^{3/2+\gep}} N ^{-2\gl q/ \Sigma_h(q)}.
 \end{split}
 \end{equation}
 The second inequality holds for $N$ sufficiently large.
 Now we set
 \begin{equation}
\delta\,  := \, \sup_{q<h } \frac{-2\gl q -\Sigma_h(q)}{\Sigma_h(q)},  
\end{equation}
and observe that
$\sup_{q<h} \left(-2\gl q -\Sigma_h(q)\right)$ is positive if and only if
$\sup_{q\in \R } \left( -2\gl q -\Sigma _h(q)\right)$ is positive.
The latter expression 
is the Legendre transform of $\Sigma_h (\cdot)$ and therefore 
it coincides with $-2\gl h+\log \M (2\gl)$, which is positive
for $h< \overh (\gl)$.
We have therefore proven
\eqref{eq:delta}. \qed

\medskip

Two remarks are in order:
\smallskip

\begin{rem}\rm
It is immediate to see that, at least in the case in which
 $\text{ess}\sup \go_1 =: \go^\star$ is finite,  
Proposition \ref{th:LBonZ} yields that C.1 cannot hold in the free endpoint case
below $\overh (\gl)$. Even more, it is in contradiction with the
possibility of having $o(\log N)$ visits to the  unfavorable
solvent, if one insists on having $\bbP(\dd \go)$--a.s. 
results. This is easily seen by considering copolymers of length $\tau _N$:
one estimates for the numerator   the 
 partition function restricted to trajectories of the walk
visiting $o(\log N)$ times the lower half plane and for
the denominator one uses \eqref{eq:delta}.
\end{rem}
\smallskip

\begin{rem}\rm
The argument in the proof of Proposition \ref{th:LBonZ}
may be repeated for the constrained case and
one obtains that, as long as $h < \omega^\star$ 
there exists $\delta>0$ such that   $N^{3/2-\delta}Z^\rc _{N ,\go}$
 does not vanish. This once again contradicts C.1: it is however a more
 evident phenomenon, since the polymer is forced in any case to
 visit all atypical $\go$--stretches when its endpoint encounters them.  
\end{rem}

\smallskip

\section*{acknowledgments} We would like to thank Thierry Bodineau 
for several enlightening discussions. We thank also Francesco Caravenna and Massimiliano Gubinelli for the interaction we had on the content of
 Section 4, and 
Michel Ledoux and C\'edric Villani
for e--mail exchanges on concentration issues. 
F.L.T. is grateful to the {\sl Laboratoire de Probabilit\'es et Mod\`eles Al\'eatoires  P6 \& 7} 
for the hospitality
and to the University of Paris 7  
for the support. This work was supported 
in part by the Swiss Science Foundation Contract No. 20-100536/1.


\begin{thebibliography}{15}  

\bibitem{cf:AZ}
S. Albeverio and X. Y. Zhou, \textit{Free energy and some sample path properties of a random walk with random potential},
  J. Statist. Phys. {\bf 83 } (1996),  573--622.
  
 \bibitem{cf:AS}
 K. S. Alexander and V. Sidoravicius,
 \textit{Pinning of polymers and interfaces by random potentials},
 preprint (2005). Available on: arXiv.org e-Print archive: math.PR/0501028 
  
 \bibitem{cf:BisdH} 
M.  Biskup and F. den Hollander, \textit{A heteropolymer near a linear interface},  Ann. Appl. Probab. {\bf 9} (1999), 668--687.

\bibitem{cf:BoGoe}
S. Bobkov and F. G\"otze,
\textit{Exponential integrability and transportation cost related to
logarithmic Sobolev inequalities}, J. Funct. Anal. {\bf 163} (1999), 1--28.

\bibitem{cf:BG}
T. Bodineau and G. Giacomin, \textit{On the localization transition   of random copolymers near selective interfaces},
J. Statist. Phys. {\bf 117} (2004), 801-818.

\bibitem{cf:BdH} E. Bolthausen and F. den Hollander,    
 \textit{Localization transition for a polymer near an interface},   
 Ann. Probab.  {\bf 25}  (1997),  1334--1366.   
   

\bibitem{cf:ch} P. Carmona and Y. Hu, {\em Universality in Sherrington-Kirkpatrick's Spin Glass Model}, 
preprint (2004). 
\\
Available online: http://www.proba.jussieu.fr/pageperso/hu/preprints.html

\bibitem{cf:DZ}  
A. Dembo and O. Zeitouni, \textit{Large deviations techniques and applications},  
Springer Verlag, New York, 1998.  
  
\bibitem{cf:DGZ}
J.--D. Deuschel, G. Giacomin and L. Zambotti,
\textit{Scaling limits of  equilibrium wetting models
in  (1+1)--dimension},
Probab. Theory Rel. Fields, 
 DOI: 10.1007/s00440-004-0401-8 (December 2004).
 Published online: http://link.springer.de/link/service/journals/00440/


\bibitem{cf:Feller1}  
W.~Feller, \textit{An introduction to probability theory and its applications}, Vol. I,  
 Third edition, John Wiley \& Sons, Inc.,   
New York--London--Sydney, 1968.  



\bibitem{cf:GHLO}
T. Garel, D. A. Huse, S. Leibler and H. Orland,
\textit{Localization transition of random chains at interfaces},
Europhys. Lett. {\bf 8} (1989), 9--13.

\bibitem{cf:G}  
G.~Giacomin, \textit{Localization phenomena in random polymer models}, preprint (2004). 
\\
Available
online: http://felix.proba.jussieu.fr/pageperso/giacomin/pub/publicat.html

\bibitem{cf:limite} F. Guerra and F.~L.~Toninelli, \textit{The Thermodynamic Limit in Mean Field
 Spin Glass Models}, Commun. Math. Phys. {\bf 230} (2002), 71--79.


\bibitem{cf:IY}
Y. Isozaki  and N. Yoshida (2001), {\it Weakly pinned random walk on the wall: pathwise descriptions of the phase transition},  Stoch. Proc. Appl. {\bf 96}, no. 2, 261--284.

\bibitem{cf:Ledoux}
M. Ledoux,
\textit{
Measure concentration, transportation cost, and functional inequalities},
Summer School on Singular Phenomena and Scaling in Mathematical Models, Bonn, 10--13 June 2003.
\\
Available online: http://www.lsp.ups-tlse.fr/Ledoux/  



\bibitem{cf:Monthus} C. Monthus, \textit{On the localization of   
random heteropolymers at the interface between two selective solvents},  
Eur. Phys. J. B {\bf 13} (2000), 111--130.   
  
  \bibitem{cf:RY}  D. Revuz and  M. Yor,
\textit{Continuous Martingales and Brownian Motion}, Springer Verlag, New York--Heidelberg 1991.


\bibitem{cf:Sinai}
Ya. G. Sinai, \textit{A 
random walk with a random potential},  Theory Probab. Appl. {\bf 38} (1993),  382--385.

\bibitem{cf:newlook} M. Talagrand, \textit{A new look at independence}, Ann. Probab. {\bf 24}  (1996), 1--34.


\bibitem{cf:talacergy} M. Talagrand, \textit{Gaussian averages, 
Bernoulli averages and Gibbs' measure}, Random Structures Algorithms {\bf 21} (2002), 197--204.

\bibitem{cf:cedric} C. Villani,
\textit{Topics in optimal transportation}, Graduate Studies in Mathematics {\bf 58}, 
American Mathematical Society, Providence, RI, 2003.





\end{thebibliography}
\end{document}